%Authors: S. Montgomery-Smith

%Title: Time decay for the bounded mean oscillation of
% solutions of the Schr\"odinger and wave equations

%Filename: montsmithschrodinger.tex

%TeX: Plain_TeX

%Length: 33485

%Received Date 4/7/97

%SubjectClass: 35J10, 35L05, 42A45, 42B30, 60J65

%Abstract: Let $u(x,t)$ be the solution of the Schr\"odinger or wave
% equation with $L_2$ initial data.  We provide counterexamples to
% plausible conjectures involving the decay in $t$ of the $\BMO$ norm
% of $u(t,\cdot)$.  The proofs make use of random methods, in particular,
% Brownian motion.
%
% Despite the applied sounding title, this really is a paper in harmonic
% analysis.

%Citation: to appear in Duke Math J.

%Special character check block
%32   space        33 ! exclam. pt.   34 " double quote  35 # sharp
%36 $ dollar       37 % percent       38 & ampersand     39 ' prime
%40 ( left paren.  41 ) rt. paren.    42 * asterisk      43 + plus
%44 , comma        45 - minus         46 . period        47 / divide
%58 : colon        59 ; semi-colon    60 < less than     61 = equal
%62 > greater than 63 ? question mark 64 @ at
%91 [ left bracket 92 \ backslash     93 ] right bracket 94 ^ caret
%95 _ underline    96 ` left single quote
%123 { left brace  124 | vertical bar 125 } right brace  126 ~ tilda

%Insert your TeX file starting here.

\magnification=\magstep1

\overfullrule = 0pt

\def\normo#1{\left\|#1\right\|}
\def\modo#1{\left|#1\right|}
\def\snormo#1{\mathopen\|#1\mathclose\|}
\def\smodo#1{\mathopen|#1\mathclose|}
\def\angleo#1{\langle#1\rangle}

\def\Re{\mathop{\hbox{\rm Re}}}
\def\sgn{\mathop{\hbox{\rm sgn}}}
\def\R{{\bf R}}
\def\BMO{{\rm BMO}}
\def\assert#1\par{\smallskip{\sl\noindent #1}\smallskip}

\centerline{\bf Time decay for the bounded mean oscillation of}
\centerline{\bf solutions of the Schr\"odinger and wave equations}

\bigskip

\centerline{S.J.~Montgomery-Smith%
\footnote*{Research funded in part by N.S.F.\ D.M.S. 9424396.}}
\centerline{University of Missouri}

\bigskip

{\narrower \noindent
Abstract: Let $u(x,t)$ be the solution of the Schr\"odinger or wave
equation with $L_2$ initial data.  We provide counterexamples to
plausible conjectures involving the decay in $t$ of the $\BMO$ norm
of $u(t,\cdot)$.  The proofs make use of random methods, in particular,
Brownian motion.

}

\bigskip

\beginsection 1. Introduction

Consider the wave equation:
$$ \eqalignno{
   \partial^2_t u(t,x) &= \Delta u(t,x) \cr
   u(0,x) &= 0 \cr
   \partial_t u(0,x) &= f ,\cr}$$
where $t \in \R$, $x \in \R^d$, 
and $f \in L_2(\R^d)$.  
Let us write once and for all
$B_t f = u(t,\cdot)$, so that $\widehat{B_t f}(\zeta) = 
{\displaystyle{\sin(t\modo\zeta) \over \modo\zeta}}\hat f(\zeta)$.  
In this paper we will be discussing the endpoints of the 
following assertion.

\assert $(B_{d,p,q})$:  if $f\in L_2(\R^d)$, then $t\mapsto B_t f$ is in
$L_q(\R,L_p(\R^d))$, and there exists a constant $c$ independent of
$f$ such that
$$ \left(\int_{-\infty}^\infty \normo{B_t f}_p^q \, dt \right)^{1/q}
   \le c \, \normo f_2 .$$

\noindent If $p,q \ge 1$, then by standard arguments, it is easy to show
that this is equivalent to its dual assertion.  Here, as in the rest of
the paper, $1/p + 1/p' = 1/q + 1/q' = 1$.

\assert $(B^*_{d,p',q'})$:  if $t\mapsto f_t$ is in 
$L_{q'}(\R,L_{p'}(\R^d))$,
then $\int_{-\infty}^\infty B_t f_t \, dt$ 
(exists almost everywhere and) is in $L_2(\R^d)$, and there
exists a constant $c$ independent of $f_t$ such that
$$ \normo{\int_{-\infty}^\infty B_t f_t\, dt}_2 \le c \,
   \left(\int_{-\infty}^\infty \normo{f_t}_{p'}^{q'} \, dt \right)^{1/q'} .$$

\noindent
In stating these assertions, we will always suppose that $d$, $p$ and 
$q$ satisfy the following conditions:
$$ \left.\eqalign{
   {d\over p} + {1\over q} &= {d\over 2} - 1 \cr
   p,q &\ge 2 \cr
   d &\ge 3 .\cr }\right\} \eqno (*_B) $$
Indeed the first condition can easily be shown to be necessary if 
$(B_{d,p,q})$
is to hold.  This follows from considering the substitution $f(x) \to f(rx)$
for any $0<r<\infty$.
(Physicists might call this a dimensional
argument.)

It is known that if conditions~$(*_B)$ hold, then assertions $(B_{d,p,q})$
and $(B^*_{d,p',q'})$ are true whenever $p<\infty$.  Results of this form
are commonly known as Strichartz inequalities, and were proven
in [St1] in the case $p = q$.  The other cases follow by the same argument.
The main application of these
results is to show existence and uniqueness results for non-linear
wave equations.  (See [Br], [Ka], [Ke] and [Ru] for applications of this and
similar results, and references to other such results.)

The only case left is $p=\infty$, $q = 2$ and $d = 3$.  The assertion
$(B_{3,\infty,2})$ and its dual $(B^*_{3,1,2})$ were shown to be false
by Klainerman and Machedon [Kl].

It was then conjectured that a weaker assertion may hold, that is,
$(B_{3,\BMO,2})$.  The
space $\BMO(\R^d)$ has long been studied by many authors.  Many results
that are not true for $L_\infty$ turn out to be true for $\BMO$.  The
space $\BMO$ is slightly larger than $L_\infty$, with a smaller norm.
There are several equivalent definitions for this space (see for example
[Se]), 
and for our
purposes it will be convenient to define $\BMO$ in terms of its dual space,
$H_1$, that is, $\BMO(\R^d)$ is the space of measurable functions $f$ from
$\R^d$ such that if $g \in H_1(\R^d)$, then $fg \in L_1(\R^d)$.  Furthermore,
$\BMO$ is equipped with the dual norm:
$$ \normo f_{\BMO} = \sup\left\{ \int_{\R^d} fg : \normo g_{H_1} \le 1
   \right\} .$$
Thus it only remains to define $H_1(\R^d)$.  Again, the literature
gives several definitions (see for example [Se]).  
We will pick a definition in terms of
maximal functions based upon the heat kernel:
$$ Mg(x) = \sup_{t>0} \int_{\R^d} {1\over (4 \pi t)^{d/2}}
   \exp(-\modo{y}^2/t) g(x-y) \, dy .$$
Then we say that $g \in H_1(\R^d)$ if $Mg \in L_1(\R^d)$, and set
$$ \normo g_{H_1} = \normo{Mg}_1 .$$
It is clear that $\modo g \le \modo{Mg}$ almost everywhere, and hence
$H_1$ is a subspace of $L_1$ with larger norm.

For definiteness, let us explicitly state the assertions involving these
norms.

\assert $(B_{3,\BMO,2})$:  if $f\in L_2(\R^3)$, then $t\mapsto B_t f$ is in
$L_2(\R,\BMO(\R^3))$, and there exists a constant $c$ independent of
$f$ such that
$$ \left(\int_{-\infty}^\infty \normo{B_t f}_{\BMO}^2 \, dt \right)^{1/2}
   \le c \, \normo f_2 .$$

\noindent This is equivalent to its dual assertion:

\assert $(B^*_{3,H_1,2})$:  if $t\mapsto f_t$ is in $L_2(\R,H_1(\R^3))$,
then $\int_{-\infty}^\infty B_t f_t \, dt$ is in $L_2(\R^3)$, and there
exists a constant $c$ independent of $f_t$ such that
$$ \normo{\int_{-\infty}^\infty B_t f_t \, dt}_2 \le c \,
   \left(\int_{-\infty}^\infty \normo{f_t}_{H_1}^2 \, dt \right)^{1/2} .$$

\noindent
The purpose of this paper is to provide a counterexample to assertion
$(B_{3,\BMO,2})$.  We will proceed by considering
the dual assertion $(B^*_{3,H_1,2})$.  We will
be using random methods, and so we do not provide an
explicit counterexample.  The tool from probability theory we shall use
is Brownian motion.  We refer the reader to [Pe] for details.  However
the only property of Brownian motion we shall use is that it is a
randomly chosen continuous function $t\mapsto b_t$ from $\R$ to $\R$
such that $b_t - b_s$ is a gaussian random variable with mean 0 and
standard deviation $\sqrt{\modo{s-t}}$.  
Throughout the paper, we will use the notation $E X$
to denote the expected value of a random variable $X$.

We shall also consider 
the Schr\"odinger equation with zero potential:
$$ \eqalignno{
   \partial _t u(t,x) &= -i \Delta u(t,x) \cr 
   u(0,x) &= f(x) , \cr }$$
where, once again, $t \in \R$, $x \in \R^d$, 
and $f \in L_2(\R^d)$.  Let us write once and for all 
$A_t f = u(t,\cdot)$, so that $\widehat{A_t f}(\zeta) = \exp(
i \modo\zeta^2) \hat f(\zeta)$.  
Then we get a similar assertion, with its dual.

\assert $(A_{d,p,q})$:  if $f\in L_2(\R^d)$, then $t\mapsto A_t f$ is in
$L_q(\R,L_p(\R^d))$, and there exists a constant $c$ independent of
$f$ such that
$$ \left(\int_{-\infty}^\infty \normo{A_t f}_p^q \, dt \right)^{1/q}
   \le c \, \normo f_2 .$$

\assert $(A^*_{d,p',q'})$:  if $t\mapsto f_t$ is in 
$L_{q'}(\R,L_{p'}(\R^d))$,
then $\int_{-\infty}^\infty A_t f_t \, dt$ is in $L_2(\R^d)$, and there
exists a constant $c$ independent of $f_t$ such that
$$ \normo{\int_{-\infty}^\infty A_t f_t \, dt}_2 \le c \,
   \left(\int_{-\infty}^\infty \normo{f_t}_{p'}^{q'} \, dt \right)^{1/q'} .$$

\noindent
In stating these assertions, we will always suppose that $p$, $q$ and $d$
satisfy the following conditions:
$$ \left.\eqalign{
   {d\over p} + {2\over q} &= {d\over 2} \cr
   p,q &\ge 2 \cr
   d &\ge 2 .\cr }\right\} \eqno (*_A) $$
Once again, the first condition can be shown to be necessary.

Rather less is known about these assertions than the corresponding
ones for the wave equation.  It is known [Gi] that if
conditions~$(*_A)$ hold, and $q > 2$, then $(A_{d,p,q})$ holds.

However the case $q = 2$ seems to be open.  This paper tackles one of these
problems, that is, the case when $q=2$, $d=2$ and $p = \infty$.  We will
demonstrate that the assertion $(A_{2,\infty,2})$ does not hold.  
(We will also include an alternative proof of this fact due to
Tony Carbery and Steve Hofmann.)
The problem
when $q = 2$ and $d \ge 3$ seems to be very difficult, and at the time
of writing is apparently unknown.

We will also deal with the assertions involving the space $\BMO$, again
showing that these are false.

\assert $(A_{2,\BMO,2})$:  if $f\in L_2(\R^2)$, then $t\mapsto A_t f$ is in
$L_2(\R,\BMO(\R^2))$, and there exists a constant $c$ independent of
$f$ such that
$$ \left(\int_{-\infty}^\infty \normo{A_t f}_{\BMO}^2 \, dt \right)^{1/2}
   \le c \, \normo f_2 .$$

\noindent This is equivalent to its dual assertion:

\assert $(A^*_{2,H_1,2})$:  if $t\mapsto f_t$ is in $L_2(\R,H_1(\R^2))$,
then $\int_{-\infty}^\infty A_t f_t \, dt$ is in $L_2(\R^2)$, and there
exists a constant $c$ independent of $f_t$ such that
$$ \normo{\int_{-\infty}^\infty A_t f_t\, dt}_2 \le c \,
   \left(\int_{-\infty}^\infty \normo{f_t}_{H_1}^2 \, dt \right)^{1/2} .$$

\noindent The author would like to acknowledge the tremendous help of many
colleagues, who explained the problems, and were a source of ideas and
inspiration.  In particular, he would like to mention 
Tony Carbery,
Loukas Grafakos,
Steve Hofmann,
Nigel Kalton,
Lev Kapitansky, and
Luis Vega.

\beginsection 2. Solutions of the Schr\"odinger equation

We will start by considering the Schr\"odinger equation, because the
techniques are simpler.  Our first result will be eclipsed by Theorem~2,
given later.  However, we will prove the following result because the method
of the proof is simpler, and illustrates the main ideas that will be used.
Because this will not be a definitive result, we will not be completely
rigorous.

Throughout this section, we will make great use of the fact that if
$\hat f(\zeta) = \exp(-\alpha \modo\zeta^2)$ ($\zeta\in\R^2$), where
$\alpha$ is a complex number with $\Re(\alpha) > 0$, then
$f(x) = (4 \pi \alpha)^{-1} \exp(-\modo x^2 /\alpha)$.  If one is considering
tempered distributions, the result remains true if $\Re(\alpha) = 0$.

\proclaim Theorem 1. Assertion $(A_{2,\infty,2})$ is not true.

In fact, what we will do is to show that assertion $(A^*_{2,1,2})$ is not
true.
Our counterexample is $f_t(x) = \alpha_t \delta(x-p_t)$, where 
$\int_{-\infty}^\infty \modo{\alpha_t}^2 \, dt = 1$, $p_t \in \R^2$ will
be chosen later,
and $\delta$ is the Dirac delta function on $\R^2$.  
Of course, $\delta$ is not a
function, and thus is not in $L_1(\R^2)$.  However, 
we will sacrifice rigor
for the sake of clarity.  A more rigorous argument may be formed by
setting $\hat \delta(\zeta) = \exp(-\modo\zeta^2)$, and following the
argument used in the proof of Theorem~2.

Note that if $t \in \R$, and $p,q \in \R^2$, then
$$ \eqalignno{
   \int_{\R^2} (A_t \delta)(x-p) \overline{\delta(x-q)} \, dx
   &=
   \int_{\R^2} \exp(it\modo\zeta^2) \exp(ip\cdot\zeta) 
   \exp(-iq\cdot\zeta) \, d\zeta \cr
   &=
   -{1\over 4 \pi i t} \exp\left({\modo{p-q}^2\over i t}\right) . \cr } $$
Then
$$  \eqalignno{
    \normo{\int_{-\infty}^\infty A_t f_t \, dt}_2^2 
    &=
    \int_{\R^2} \int_{-\infty}^\infty A_t f_t(x) \, dt
    \int_{-\infty}^\infty \overline{A_s f_s(x)} \, ds \, dx \cr
    &=
    \int_{-\infty}^\infty \int_{-\infty}^\infty
    \int_{\R^2} (A_{t-s} f_t)(x) \overline{f_s(x)} \, dx \, ds \, dt \cr
    &=
    -\int_{-\infty}^\infty \int_{-\infty}^\infty
    { \alpha_t \bar \alpha_s  \over 4 \pi i (t-s)}
    \exp\left({\modo{p_t - p_s}^2 \over i(t-s)}\right) \, ds \, dt .\cr } $$
Now we let $p_t = (\sqrt\theta b_t,0)$, where $\theta$ is chosen so that 
if $\gamma$ is a gaussian random variable with mean 0 and standard deviation
$1$, then $a = E(\sin(\theta\gamma^2) \ne 0$.  
Thus 
$$ E \sin\left({\modo{p_t - p_s}^2 \over (t-s)}\right)
   = E \sin(\sgn(t-s) \theta \gamma^2)
   = a \sgn(t-s) .$$
Then, if $\alpha_t$ is real,
$$ \eqalignno{
   -\Re E\normo{\int_{-\infty}^\infty A_t f_t \, dt}_2^2 
   &=
   \int_{-\infty}^\infty \int_{-\infty}^\infty
    { \alpha_t \alpha_s  \over 4 \pi i (t-s)}
    E\sin\left({\modo{p_t - p_s}^2 \over (t-s)}\right) \, ds \, dt \cr
    &=
    a \, \int_{-\infty}^\infty \int_{-\infty}^\infty
    { \alpha_t \alpha_s  \over 4 \pi i \modo{t-s}}
    \, ds \, dt .\cr }$$
Taking $\alpha_t = 1$ if $\modo{t} \le 1$, and $0$ otherwise, we
obtain an unbounded integral, and hence the desired counterexample.

\bigskip

We will also give a different proof of this last result due
to Tony Carbery and Steve Hofmann.  This proof was found shortly after
the one just given, and is reproduced here with their permission.  
Again, we are sacrificing some rigor to obtain
clarity.
The starting point is the same, except that we shall
take $f_t(x) = \delta(x-\sqrt r p_t)$ for $t\in I$, and $0$ otherwise, 
where $I$ is any 
interval.  We shall suppose that
$p_t$ is any path in $\mathord{\hbox{\rm Lip}}_{1/2}$, with $\normo p_
{\mathord{\hbox{\sevenrm Lip}}_{1/2}} = 1$,
and that $r \in \R$.  Then it is sufficient to find
a counterexample to the following assertion.  For every $r\in\R$ and
interval $I$ we have
$$ \modo{\int_I \int_I {1\over 4\pi i(t-s)} 
   \exp \left({r \modo{p_t-p_s}^2
   \over i(t-s)}\right) 
   \, ds \, dt} \le \modo{I} .$$
Now we apply an idea from [Co].
Pick $F\in C^\infty_0(\R)$ such that $F(\rho) = \rho$ for $\modo \rho
\le 1$.  Then
$$ {\modo{p_t-p_s}^2\over t-s} =
   \int_{-\infty}^\infty \hat F(-r) 
   \exp \left({r \modo{p_t-p_s}^2
   \over i(t-s)}\right) \, dr ,$$
and so combining the last two displayed equations, and rearranging
the integrals, 
we see that
$$ \modo{\int_I \int_I {\modo{p_t-p_s}^2 \over 4\pi (t-s)^2} 
   \, ds \, dt} \le \modo{I} \int_{-\infty}^\infty \smodo{\hat F(r)} \, dr
   \le C \modo I ,$$
where $C$ is some universal constant.  But by a result of Strichartz [St2], 
we have that if $D^{1/2} p \in \BMO$, then
$$ \sup_I {1\over \modo I} \int_I \int_I
   {\modo{p_t-p_s}^2 \over (t-s)^2} \, ds \, dt
   \approx
   \snormo{ D^{1/2} p}_{\BMO} .$$
Here $\widehat{D^{1/2} p}(\tau) = \modo \tau^{1/2} \hat p(\tau)$.
Hence we have our counterexample by picking $p$ with
$\normo p_{\mathord{\hbox{\sevenrm Lip}}_{1/2}} = 1$, and 
$\snormo{D^{1/2} p}_{\BMO}$ arbitrarily large.

\bigskip

Now we will improve this result.  In the arguments that
follow, the reader may feel uncomfortable with the cavalier
and implicit use of Fubini's Theorem.  The use of Fubini's Theorem
requires all the integrals to be absolutely convergent, and this is not
the case, as it is
exactly the opposite that we are trying to show.  For this reason, the
proofs should really be seen as proofs by contradiction, that is, the
reader should suppose initially that the assertions stated in the theorems
are true.  

\proclaim Theorem 2.  Assertion $(A_{2,\BMO,2})$ is not true.

Once again, we provide a counterexample to the dual 
assertion $(A^*_{2,H_1,2})$.
In this case, our counterexample will be $f_t(x) = \alpha_t g(x-p_t)$, where
$\int_{-\infty}^\infty \modo{\alpha_t}^2 \, dt = 1$, $p_t \in \R^2$ will
be chosen later,
and $\hat g(\zeta) = \modo{\zeta}^2 \exp(-\modo\zeta^2)$.  

Let us first show that $g \in H_1(\R^2)$.  
Using Fourier transforms, we see that
$$ \eqalignno{
   Mg(x) 
         &= 
         \sup_{t>0} {1\over 4\pi t} \int_{\R^2} \exp(-\modo{x-y}^2/t)
         g(y) \, dy \cr
         &= 
         \sup_{t>0} {4(1+t) - 2 \modo x^2 \over 4\pi(1+t)^3} 
           \exp\left(-{\modo x^2\over 1+t}\right) \cr
         &= 
         O({1\over 1+\modo x^2}) ,\cr} $$
which is in $L_1(\R^2)$.

Note also
that for $t\in\R$ and $p,q\in\R^2$
$$ \eqalignno{
   \int_{\R^2} & (A_t g)(x-p) \overline{g(x-q)} \, dx \cr
   &=
   \int_{\R^2} \modo{\zeta}^4 \exp((-2 + i t)\modo{\zeta}^2) \exp(\zeta
   \cdot(q-p)) \, d\zeta \cr
   &=
   {1\over 4\pi (2-it)} {\Delta^2} 
   \exp\left({\modo x^2 \over -2 + it}\right) \Big|_{x = q-p} \cr
   &=
   {1\over 4\pi (2-it)^3} P(\modo{q-p}^2/(-2+it))
   \exp\left({\modo{q-p}^2 \over -2 + it}\right) \cr
   &=
   {1\over 4\pi (2-it)^3} P(\modo{q-p}^2/(-2+it))
   \exp\left(-{2\modo{q-p}^2 \over 4 + t^2}
   - {it\modo{q-p}^2 \over 4 + t^2}\right), \cr }$$
where $P(t) = 16t^2 - 64t + 32$. 

Hence
$$ \normo{\int_{-\infty}^\infty A_t f_t \, dt}_2^2 
   =
   \int_{-\infty}^\infty \int_{-\infty}^\infty
   \alpha_t \bar \alpha_s K(t,s) \, ds \, dt ,$$
where
$$ \eqalignno{
   K(t,s) &= 
   {1 \over 4 \pi (2-i(t-s))^3}
   P(\modo{p_t-p_s}^2/(-2+i(t-s))) \times \cr
   & \exp\left(-{2\modo{p_t-p_s}^2 \over 4 + (t-s)^2}
   - {i(t-s)\modo{p_t-p_s}^2 \over 4 + (t-s)^2}\right). \cr }$$
Now we will choose $p_t = (t,\sqrt\theta b_t)$, where $b_t$ and $\theta$ 
are defined as in the previous section.  In order to demonstrate
that we have a counterexample, we need to show that 
$E(K(t,s))$ is not the kernel of  bounded operator from $L_2(\R) \to
L_2(\R)$.

We shall be interested in the behavior of $E(K(t,s))$ as
$t-s \to \pm \infty$.  Let $b_t - b_s = \sqrt{\modo{t-s}} \gamma$, where
$\gamma$ is a gaussian random variable with mean $0$ and standard 
deviation $1$.  Denote $\angleo{t} = 1 + \modo t$.

First see that
$$ P(\modo{p_t-p_s}^2/(-2+i(t-s))) = -16(t-s)^2 + O(\gamma^4 \angleo{t-s}).$$
Also, there exists a constant $c_1$ such that
$$ \exp\left(-{2\modo{p_t-p_s}^2 \over 4 + (t-s)^2}\right) 
   =
   \exp(-(2 + O(\gamma^2\angleo{t-s}^{-1})) 
   =
   \exp(-2) + O(\exp(c_1 \gamma^2\angleo{t-s}^{-1})-1).$$
Furthermore,
$$ {(t-s) \modo{p_t-p_s}^2 \over 4 + (t-s)^2}
   = 
   (t-s) + \sgn(t-s)\theta\gamma^2 + O(\gamma^2\angleo{t-s}^{-1}) ,$$
and hence there is a constant $c_2$ such that
$$ \exp\left(-i {(t-s) \modo{p_t-p_s}^2 \over 4 + (t-s)^2}\right)
   =
   \exp(i(s-t + \sgn(s-t)\theta\gamma^2))
   + O(\exp(c_2 \gamma^2\angleo{t-s}^{-1})-1) .$$
Therefore, for some constant $c_3$,
$$ \eqalignno{
   K(t,s) 
   &= 
   {-16 \exp(-2) (t-s)^2 \over 4 \pi (-2 + i(t-s))^3}
   \exp(i(s-t)) \exp(i\sgn(s-t) \theta\gamma^2) \cr
   &+
   O\bigl(\gamma^4\angleo{t-s}^{-3} +
   \angleo{t-s}^{-2} (\exp(c_3 \gamma^2\angleo{t-s}^{-1})-1)\bigr) .\cr }$$
Hence
$$ E(K(t,s)) =
   {-16 \exp(-2) (t-s)^2 \over 4\pi (2-i(t-s))^3}
   \exp(i(s-t)) (a_1 + i a_2 \sgn(s-t))
   + O(\angleo{t-s}^{-3}) ,$$
where $a_1 + i a_2 = E(\exp(i \theta \gamma^2))$.  
By considering the examples $\alpha_t = \exp(it)/\sqrt{N}$ if $\modo{t} \le
N$, and $0$ otherwise, and letting $N \to \infty$, it may be readily
seen that
$E(K(t,s))$ is not the kernel of a bounded
map from $L_2(\R) \to L_2(\R)$, and the desired counterexample has been
obtained.

\bigskip

We will now present a second proof of the same result.  
This proof goes via the Fourier transform.  Although this proof is
less intuitive, it also has less technical difficulties.
To recap, it is 
sufficient to find a counterexample to the following statement:
there is a constant $c$ such that
if $\int_{-\infty}^\infty \modo{\alpha_t}^2 \, dt \le 1$, and if
$g \in H_1(\R^2)$, then for any path $p_t \in \R^2$, we have that
$$ \int_{-\infty}^\infty \int_{-\infty}^\infty 
   K(s,t) \alpha_t \bar \alpha_s \, ds \, dt \le c,$$
where
$$ \eqalignno{
   K(s,t) 
   &=
   \int_{\R^2} (A_{t-s} g)(x-p_t) \overline{g(x-p_s)} \, dx \, ds \, dt \cr
   &=
   \int_{\R^2} \exp(i (t-s) \modo\zeta^2 + i(p_s-p_t)\cdot\zeta) 
   \modo{\hat g(\zeta)}^2
   \, d\zeta .\cr}$$
That is, we are asking whether $K(s,t)$ is the kernel
of a bounded operator from $L_2(\R)$ to $L_2(\R)$.

We take $p_t = (t,b_t)$, where $b_t$ is Brownian motion.  It
is sufficient to show that $E(K(s,t))$ is not the kernel of a bounded
operator on $L_2(\R^2)$.   
Note that if $\gamma$ is a gaussian random variable with mean $1$ and
standard deviation $0$, then $E(\exp(i\gamma)) = \exp(-t^2/2)$.  Hence
$ E(K(s,t)) = L(t-s)$, where
$$ L(t) =
   \int_{\R^2} \exp(i t \modo\zeta^2 - i t\zeta_1 - \modo t \zeta_2^2/2)
   \modo{\hat g(\zeta)}^2 \, d\zeta .$$
But $L(t-s)$ fails to be the kernel of a bounded operator on $L_2(\R)$ if
and only if $\hat L(\omega)$ fails to be bounded almost everywhere.  
Furthermore
$$ \hat L(\omega) =
   \int_{\R^2} k(\omega,\zeta) \modo{g(\zeta)}^2 \, d\zeta ,$$
where $k(\cdot,\zeta) = k_1(\cdot,\zeta)*k_2(\cdot,\zeta)$ with
$$ k_1(\omega,\zeta) = \delta(\omega + \modo\zeta^2 - \zeta_1) $$
and
$$ k_2(\omega,\zeta) = {4\zeta_2^2 \over 4\omega^2 + \zeta_2^4} ,$$
that is,
$$ k(\omega,\zeta) = 
   {4\zeta_2^2 \over 4(\omega + \modo\zeta^2 - \zeta_1)^2 + \zeta_2^4} .$$
Clearly $k$ enjoys enough continuity properties so that it is sufficient
to show that $\hat L(\omega)$ is unbounded for $\omega = 1/4$.  So
$$ k(1/4,\zeta) =
   {4\zeta_2^2 \over 4((\zeta_1-1/2)^2 + \zeta_2^2)^2 + \zeta_2^4}
   \ge
   {\zeta_2^2 \over 2\modo{\zeta - (1/2,0)}^4 } .$$
It is clear that $k(1/4,\zeta)$ has an $L_1(\R^2)$ singularity at
$(1/2,0)$, and hence taking $\hat g(\zeta) = \modo\zeta^2 
\exp(-\modo\zeta^2)$, we see that $\hat L(1/4)$ is unbounded.

\beginsection 3. Solutions of the wave equation

This section is devoted to the following result.

\proclaim Theorem 3.  Assertion $(B_{3,\BMO,2})$ is not true.

We will show that assertion $(B^*_{3,H_1,2})$ is not
true.  The methods will be essentially the same as in the previous
proofs, but the details will
be more difficult, and so we will break it into steps.

To start, let us define an operator on a dense subspace of
$L_2(\R^3)$ for each
$t\in \R$ given by
$$ \widehat{C_t f}(\zeta) =  {\cos(t\modo\zeta) \over \modo\zeta^2 }
   \hat f(\zeta) .$$
Let us also set
$$ K_t(x) = {1\over 4 \pi \modo x} 
   I_{\modo x \ge \modo t} \quad (x\in\R^3) .$$

\proclaim Lemma 4.  If $f \in L_1(\R^3) \cap L_2(\R^3)$, and
$t \ne 0$, then
$$ C_t f(x) = \int_{\R^3} K_t(x-y) f(y) \, dy ,$$
and the operator norm of $C_t$ from $L_1(\R^3)$ to $L_1(\R^3)$ is
$(4\pi\modo t)^{-1}$.

To show this when $f$ is $C^\infty$ with compact support, 
it is sufficient to show that
$\hat K_t = \cos(t\modo\zeta) / \modo\zeta^2 $ (as tempered distributions).  
Then by Young's convolution formula, as an operator
from $L_1(\R^3)$ to $L_1(\R^3)$, the operator norm of $C_t$ is given
by
$$ \normo{C_t} = \normo{K_t}_\infty = {1\over 4\pi\modo t} .$$
It is
clear that $\hat K_t(\zeta)$ depends only upon $t$ and
$\modo\zeta$.  So without loss of generality, we may suppose that
$\zeta = (\zeta_1,0,0)$, where $\zeta_1 = \modo\zeta$.  In performing
the following integral, we will use the following change of variables:
$x_1 = u$, $x_2 = \sqrt{v^2-u^2}\cos(\theta)$, and $x_3 = \sqrt{v^2-u^2}
\sin(\theta)$.  Thus the Jacobian $\partial(x_1,x_2,x_3)/\partial(u,v,\theta)
= v = \modo x$.  Then
$$ \eqalignno{
   \hat K_t(\zeta)
   &=
   \int_{\modo x \ge \modo t} 
   {\exp(-i x\cdot\zeta) \over 4 \pi \modo x}
   \, dx \cr
   &=
   {1\over 4\pi} \int_{v=\modo t}^\infty \int_{u=-v}^v \int_0^{2\pi}
   \cos(u \modo \zeta) 
   \, d\theta \, dv \, du \cr
   &=
   {\cos(t \modo\zeta) \over \modo\zeta^2 } . \cr } $$
In the last line we have used the assertion that $\lim_{R \to \infty}
\cos(R\modo\zeta) = 0$, which is true in the space of tempered distributions.
   
\proclaim Corollary 5.  If $f \in L_2(\R^3)$ is bounded with
compact support, and $s,t\in\R$, then
$$ B_s B_t f(x) 
   = {C_{s+t}f(x) - C_{s-t}f(x)\over 2}
   = \int_{\R^3} {K_{s+t}(x-y) - K_{s-t}(x-y) \over 2}
     f(y) \, dy .$$

This follows because $\sin(t\modo\zeta)\sin(s\modo\zeta)
= {1\over 2} (\cos((s+t)\modo\zeta) - \cos((s-t)\modo\zeta))$.

\proclaim Lemma 6. There is a universal constant $c$ such that 
if $t\mapsto f_t$ is in $L_2(\R,L_1(\R^3))$, then
$$ \int_0^\infty \int_0^\infty \int_{\R^3} \int_{\R^3}
   K_{s+t}(x-y) f_t(x) f_s(y)
   \, dy \, dx \, dt \, ds
   \le c \, 
   \int_0^\infty \normo{f_t}_1^2 \, dt .$$

The proof of this result depends upon the boundedness of the Hardy
operators defined on $L_2([0,\infty))$:
$$ \eqalignno{
   H \alpha(t) &= {1\over t} \int_0^t \alpha(s) \, ds \cr
   H^*\alpha(t) &= \int_t^\infty {\alpha(s)\over s} \, ds . \cr }$$
If $\alpha \in L_2([0,\infty))$, then both $H\alpha$ and $H^*\alpha$
are in $L_2([0,\infty))$, with
$\normo{H\alpha}_2 , \normo{H^*\alpha}_2 \le 2 \normo{\alpha}_2$.
The result for $H$ may be found
in [Ha], and the result for $H^*$ follows because $H^*$ is the 
adjoint operator to $H$.

Next, by Lemma~4, if $s,t > 0$
$$ \int_{\R^3} \int_{\R^3} K_{s+t}(x-y) f_t(x) f_s(y) \, dy \, dx
   \le 
   {1\over 4\pi(s+t)} \normo{f_t}_1 \normo{f_s}_1 .$$
Thus, setting $\alpha(t) = \normo{f_t}_1$, we see that
$$ \eqalignno{
   \int_0^\infty &\int_0^\infty \int_{\R^3} \int_{\R^3}
   K_{s+t}(x-y) f_t(x) f_s(y)
   \, dy \, dx \, dt \, ds \cr
   &\le 
   {1\over 4\pi} \int_0^\infty (H\alpha(t) + H^*\alpha(t)) \alpha(t) \, dt
   \cr
   &\le 
   {1\over \pi} \normo{\alpha}_2^2 .\cr } $$

\bigskip

Now we will consider the following assertion.

\assert $(C)$:  There is a universal constant $c$ such that
if $t\mapsto f_t$ is in $L_2(\R,H_1(\R^3))$, then
$$ \int_{-\infty}^\infty \int_{-\infty}^\infty \int_{\R^3} \int_{\R^3}
   K_{s-t}(x-y) f_t(x) f_s(y) \, dy \, dx \, dt \, ds
   \le 
   c \, \int_{-\infty}^\infty \normo{f_t}_{H_1}^2 \, dt . $$

\noindent
Let us first demonstrate that the failure of assertion $(C)$ implies the
failure of assertion $(B_{3,H_1,2})$.  Let us suppose that we have a
sequence $t\mapsto f_t$ in $L_2(\R,H_1(\R^3))$ such that 
$$ I = \int_{-\infty}^\infty \int_{-\infty}^\infty \int_{\R^3} \int_{\R^3}
   K_{s-t}(x-y) f_t(x) f_s(y) \, dy \, dx \, dt \, ds $$
is unbounded.  We may split the integral up into four pieces:
$$ \eqalignno{
   I &= I_1 + I_2 + I_3 + I_4 \cr
     &= \int_0^\infty \int_0^\infty \dots \, dt \, ds
     + \int_{-\infty}^0 \int_0^\infty \dots \, dt \, ds
     + \int_0^\infty \int_{-\infty}^0 \dots \, dt \, ds
     + \int_{-\infty}^0 \int_{-\infty}^0 \dots \, dt \, ds .\cr}$$
By Lemma~6, we know that $I_2$ and $I_3$ are bounded.  Therefore one of
$I_1$ or $I_4$ is unbounded, and without loss
of generality we may suppose that $I_1$ is unbounded.  So without loss
of generality, we may suppose that $f_t = 0$ if $t < 0$.

Now, multiplying out, and applying Corollary~5, we see that
$$ \eqalignno{
   \normo{ \int_{-\infty}^\infty B_t f_t \, dt}_2^2
   &=
   \normo{ \int_0^\infty B_t f_t \, dt}_2^2 \cr
   &=
   \int_0^\infty \int_0^\infty
   \int_{\R^3} \int_{\R^3}
   {K_{s+t}(x-y) - K_{s-t}(x-y)\over 2} f_t(x) f_s(y)
   \, dy \, dx \, dt \, ds . \cr}$$
Again, by Lemma~6, the last integral differs from $I_1/2$ 
by a bounded amount,
and we have produced the desired counterexample.

\bigskip

All that remains to be shown is the following.

\proclaim Lemma 7.  Assertion $(C)$ is false.

Our counterexample will be $f_t(x) = \alpha_t g(x-p_t)$, where 
$\int_0^\infty \modo{\alpha_t}^2 \, dt = 1$, and $p_t = (t,b_t,b_t')$.
The function $g \in H_1(\R^3)$ will be selected later.  
Here $b_t$ and $b'_t$ are two independent Brownian motions.
We will show that
$$ J = E \left(
   \int_{-\infty}^\infty \int_{-\infty}^\infty \int_{\R^3} \int_{\R^3}
   K_{s-t}(x-y) \alpha_t \bar\alpha_s g(x-p_t) \overline{g(x-p_s)}
   \, dy \, dx \, dt \, ds \right) $$
cannot be universally bounded.  This quantity is more easily
computed using the Fourier transform:
$$ J = E \left(
   \int_{-\infty}^\infty \int_{-\infty}^\infty
   \alpha_t \bar\alpha_s 
   \int_{\R^3} {\cos((t-s)\modo\zeta) \over \modo\zeta^2} 
   \exp(i(p_s-p_t)\cdot \zeta) \modo{\hat g(\zeta)}^2
   \, d\zeta \, dt \, ds \right) .$$
Notice that
$$ \eqalignno{
   E (\exp(i(p_s-p_t))
   &=
   \exp(i(s-t)\zeta_1) E(\exp(i(b_s-b_t)\zeta_2)+i(b_s'-b_t')\zeta_3) \cr
   &=
   \exp(i(s-t)\zeta_1) \exp(-\modo{s-t} 
   (\zeta_2^2 + \zeta_3^2) / 2) .\cr}$$
Hence
$$ J = \int_{-\infty}^\infty \int_{-\infty}^\infty
   L(s-t) \alpha_t \bar\alpha_s \, dt \, ds ,$$
where
$$ L(t) = \int_{\R^3} {\cos(t\modo\zeta) \over \modo\zeta^2}
   \exp(it\zeta_1) \exp(-\modo t (\zeta_2^2 + \zeta_3^2) / 2) 
   \modo{\hat g(\zeta)}^2 \, d\zeta.$$
Saying that $J$ is bounded for all $\alpha_t$ in $L_2(\R)$ is equivalent
to saying that convolving with $L$ gives a bounded operator on $L_2(\R)$.
But this is equivalent to $\hat L$ being in $L_\infty(\R)$.
Now
$$ \hat L(\omega) = 
   \int_{\R^3}
   k(\omega,\zeta) \modo{\hat g(\zeta)}^2/\modo\zeta^2 \, d\zeta ,$$
where $k(\cdot,\zeta) = h_1(\cdot,\zeta)*h_2(\cdot,\zeta)*h_3(\cdot,\zeta)$
with
$$ \eqalignno{
   h_1(\omega,\zeta)
   &=
   {\delta(\omega+\modo\zeta) + \delta(\omega-\modo\zeta) \over 2} \cr
   h_2(\omega,\zeta)
   &=
   \delta(\omega-\zeta_1) \cr
   h_3(\omega,\zeta)
   &=
   {4(\zeta_2^2+\zeta_3^2) \over 4\omega^2 + (\zeta_2^2+\zeta_3^2)^2 } 
   .\cr}$$
Thus
$$ k(\omega,\zeta) = {k_1(\omega,\zeta) + k_2(\omega,\zeta)\over 2} ,$$
where
$$ k_1(\omega,\zeta) =
   { 4(\zeta_2^2 + \zeta_3^2)
   \over 4(\omega - \zeta_1 + \modo\zeta)^2 + (\zeta_2^2 + \zeta_3^2)^2} ,$$
and
$$ k_2(\omega,\zeta) =
   { 4(\zeta_2^2 + \zeta_3^2)
   \over 4(\omega - \zeta_1 - \modo\zeta)^2 + (\zeta_2^2 + \zeta_3^2)^2} .$$
Since all the expressions involved are
positive, we will have found a counterexample if we can show that
$$ \int_U k_1(0,\zeta) 
   \modo{\hat g(\zeta)}^2/\modo{\zeta^2} \, d\zeta = \infty ,$$
where $U$ is any subset of $\R^3$.
We will take 
$$ U = \{ \zeta : \sqrt{\zeta_2^2+\zeta_3^2} \le \zeta_1 \le 1 \} ,$$
and 
$$ g(x) = \cases{
   1 & if $0<x_1\le 1$ and $-1\le x_2,x_3 \le 1$\cr
   -1 & if $-1\le x_1<0$ and $-1\le x_2,x_3 \le 1$\cr
   0 & otherwise,\cr } $$
so that
$$ \hat g(\zeta) = {8 (1-\cos(\zeta_1)) \sin(\zeta_2) \sin(\zeta_3)
   \over \zeta_1 \zeta_2 \zeta_3 } .$$
In this case, we see that $\modo{\hat g(\zeta)}^2 / \modo{\zeta}^2$ is
bounded from below by a positive number on $U$.  
Therefore, we need to show that
$$ \int_U k_1(0,\zeta) \, d\zeta = \infty .$$
To do this, let us compute the integral using cylindrical coordinates,
that is, we will write $z = \zeta_1$, and $r = \sqrt{\zeta_2^2+\zeta_3^2}$.
Then we see that the last integral is
$$ \int_{z=0}^1 \int_{r=0}^z 
   {4r^2 \over 4(\sqrt{z^2+r^2} - z)^2 + r^4} 2\pi r \, dr \, dz .$$
But in $U$,
$$ \sqrt{z^2+r^2} - z \le r^2/z ,$$
and hence
$$ {4r^2 \over 4(\sqrt{z^2+r^2} - z)^2 + r^4}
   \ge
   {4z^2 \over r^2 (4+z^2) } .$$
Hence the integral in question becomes bounded below by
$$ \int_{z=0}^1 \int_{r=0}^z {8 \pi z^2 \over r(4+z^2)}
   \, dr \, dz ,$$
and this is easily seen to be infinite.

\beginsection 4. References

There is an extensive literature on the positive results, and we do
not desire to give a complete list.  Instead, we refer the reader to
the papers chosen principally because they are recent, and
provide many references on past work.

\bigskip

\item{Br}Ph.~Brenner, On $L_p$--$L_{p'}$ estimates for the wave equation,
{\sl Math Z., {\bf 145}, (1975), 251--254.}
\item{Co}R.R.~Coifman, G.~David and Y.~Meyer, La solution des conjecture
de Calder\'on,
{\sl Advances in Math., {\bf 48}, (1983), 144--148.}
\item{Gi}J.~Ginibre and G.~Velo, Scattering theory in the energy space
for a class of non-linear Schr\"odinger equations, {\sl J.\ Math.\
Pure Appl., {\bf 64}, (1985), 363--401.}
\item{Ha}G.H.~Hardy, J.E.~Littlewood and G.~P\'olya, {\sl Inequalities},
(2nd.\ Ed.),
Cambridge University Press (1952).
\item{Ka}L.~Kapitanski, Global and unique weak solutions of nonlinear
wave equations, {\sl Math.\ Research Papers, {\bf 1}, (1994), 211--223.}
\item{Ke}C.E.~Kenig, G.~Ponce and L.~Vega, On the IVP for nonlinear
Schr\"odinger equations, {\sl Contemp.\ Math., {\bf 180}, (1995), 353--367.}
\item{Kl}S.~Klainerman and M.~Machedon, Space time estimates for null
forms and the local existence theorem,
{\sl Communications of Pure and Applied Math., {\bf 46}, (1993), 1221--1268.}
\item{Pe}K.E.~Peterson, {\sl Brownian Motion, Hardy Spaces and Bounded
Mean Oscillation}, Cambridge University Press, 1977.
\item{Ru}A.~Ruiz and L.~Vega, Local regularity of solutions to wave
equations with time-dependent potentials, {\sl Duke Math.\ J., {\bf 76},
(1994), 913--940.}
\item{Se}E.M.~Stein, {\sl Harmonic Analysis}, Princeton University Press,
1993.
\item{St1}R.S.~Strichartz, A priori estimates for the wave equation and some 
applications, {\sl J.\ Funct.\ Anal., {\bf 5}, (1970), 218--235.}
\item{St2}R.S.~Strichartz, BMO and Sobolev spaces,
{\sl Indiana U.\ Math.\ J., {\bf 29}, (1980), 539-558.}

\bigskip

S.J.~Montgomery-Smith

Mathematics Department

University of Missouri

Columbia, MO 65211

E-mail: stephen@math.missouri.edu

\bye